\documentclass[12pt]{amsart}
\usepackage{xy,euscript,somedefs,tracefnt,latexsym,,rawfonts,epsfig,amsfonts,amssymb,latexsym,enumerate,amscd}
\input xy
\xyoption{all}
\def\cal{\mathcal}
\def\Bbb{\mathbb}
\def\g{\gamma}

\def\G{\Gamma}

\def\r{\rangle}
\def\l{\langle}
\def\t{\times}

\newtheorem{thm}{Theorem}[section]

\newtheorem{lemma}{Lemma}[section]
\newtheorem{cor}{Corollary}[section]

\newtheorem{defn}{Definition}[section]

\newtheorem{rem}{Remark}[section]
\numberwithin{equation}{section}
\begin{document}
\thanks{To appear in {\it Archiv der Mathematik}}
\date{February 18, 2011.}
\title[Surgery groups of braid groups]{Surgery groups of the 
fundamental groups of hyperplane arrangement complements}
\author[S. Roushon]{S. Roushon}
\address{School of Mathematics, Tata Institute, 
Homi Bhabha Road, Mumbai-400005, India.}
\begin{abstract} Using a recent result of Bartels and 
L\"{u}ck (\cite{BL}) we deduce that the Farrell-Jones Fibered Isomorphism 
conjecture in $L^{\l -\infty \r}$-theory is true for any group which contains  
a finite index strongly poly-free normal subgroup, in particular, for 
the Artin full braid groups. As a consequence we explicitly compute the  
surgery groups of the Artin pure braid groups. This is obtained as a corollary  
to a computation of the surgery groups of a more general class of groups, 
namely for the fundamental group of the 
complement of any fiber-type hyperplane arrangement in ${\Bbb C}^n$.
\end{abstract} 



\maketitle


\section{Introduction}
The purpose of this short note is to compute explicitly the 
surgery ($L$-)groups of the Artin pure braid groups ($PB_n$). 
This computation requires the solutions of two other 
problems. Firstly, one has to compute the lower algebraic $K$-theory 
of the group and secondly, to show that the classical assembly map in $L$-theory 
is an isomorphism. This gives an interpretation of the 
surgery groups in terms of a generalized homology theory.

For $PB_n$ we had already computed 
the lower algebraic $K$-theory in \cite{AFR}. In \cite{FR} 
we computed it for any subgroup of the Artin full braid group ($B_n$), 
in particular for any subgroup of $PB_n$.  
Here we show that the classical assembly map in $L$-theory is  
an isomorphism for any subgroup of $B_n$. The main ingredients behind the proof 
is the $K$-theoretic vanishing result [\cite{FR}, theorem 1.1] and a recent 
result of Bartels and L\"{u}ck ([\cite{BL}, theorem B]). The later result is 
used to show that the $L^{\l -\infty \r}$-theory Fibered Isomorphism conjecture of 
Farrell and Jones ([\cite{FJ}, \S1.7]) is true for any subgroup of $B_n$.  
Finally, using the stable homotopy type 
of the corresponding Eilenberg-Maclane space of $PB_n$ 
from \cite{S} we do the computation of the 
surgery groups.

Let us first say a few words about the group $PB_n$ (and $B_n$)    
before we state the computation of its surgery groups. 
The {\it Artin full braid group} $B_n$ is generated by the 
symbols $\sigma_1,\ldots , \sigma_n$  
with respect to the relations $\sigma_i\sigma_j=\sigma_j\sigma_i$ for 
$|i-j|\geq 2$ and $\sigma_i\sigma_{i-1}\sigma_i=\sigma_{i-1}\sigma_i\sigma_{i-1}$ 
for $i\leq n$. See \cite{EA} for the original source of this group. 
The map $q:B_n\to S_{n+1}$ sending the generator $\sigma_i$ to 
the transposition $(i, i+1)$ defines a homomorphism onto the 
symmetric group $S_{n+1}$ on $(n+1)$-symbols and the kernel of this 
homomorphism is defined as the {\it Artin pure braid group} $PB_n$. 
See \cite{Bir} and \cite{KT} for some more information on braid groups. 

In this paper we need a topological 
interpretation of $PB_n$ which we describe below.
 
Let ${\cal H}_n$ be the hyperplane arrangement  
$H_{ij}=\{(x_0,x_1\ldots , x_n)\in{\Bbb C}^{n+1}\ |\ x_i=x_j\}$ for 
$i,j=0,1,\ldots , n$ and $i\neq j$ in the $(n+1)$-dimensional 
complex space ${\Bbb C}^{n+1}$.
The fundamental group of the complement 
${\Bbb C}^{n+1}-\cup_{i\neq j}H_{ij}$ 
is isomorphic to $PB_n$. Note that the group $S_{n+1}$ 
acts freely on ${\Bbb C}^{n+1}-\cup_{i\neq j}H_{ij}$ by 
permuting coordinates. The fundamental group of the quotient space 
$({\Bbb C}^{n+1}-\cup_{i\neq j}H_{ij})/S_{n+1}$ is isomorphic to $B_n$. 
Therefore, there is an exact sequence of the following type. 

$$1\to PB_n\to B_n\to S_{n+1}\to 1.$$

Furthermore, the homomorphism $B_n\to S_{n+1}$ in the above exact 
sequence coincides with $q$ we defined above. 
For this interpretation of the braid groups see \cite{FN}.

\begin{cor}\label{maincor} 
For all $n\geq 1$ the surgery groups of the Artin pure braid group 
$PB_n$ are computed as follows.

$$
L_i(PB_n) =\begin{cases}{\Bbb Z}&\text{if
$i\equiv$ 0 mod 4}\\
{\Bbb Z}^{\frac{n(n+1)}{2}} &\text{if $i\equiv$ 1 mod 4}\\
{\Bbb Z}_2& \text{if $i\equiv$ 2 mod 4}\\
{\Bbb Z}_2^{\frac{n(n+1)}{2}}& \text{if $i\equiv$ 3 mod 4}.
\end{cases}$$
\end{cor}

\begin{proof} This is an immediate corollary of 
Theorem \ref{mainthm2} since the arrangement ${\cal H}_n$ is fiber-type and 
there are $\frac{n(n+1)}{2}$ 
hyperplanes in ${\cal H}_n$.
\end{proof}

We recall here that there are surgery groups for different kinds of 
surgery problems and 
they appear in the literature with the notations 
$L^*_i(-)$, where $*=h,s, \l -\infty \r$ or
$\l i\r$ for $i\leq 0$. But all of them are naturally 
isomorphic for torsion-free groups 
$G$ if the Whitehead group $Wh(G)$, the reduced 
projective class group $\tilde K_0({\Bbb Z}G)$ and the  
negative $K$-groups $K_{-i}({\Bbb Z}G)$ for $i\geq 1$ 
vanish. See [\cite{LR}, remark 1.21 and 
proposition 1.23]. 
Therefore, we use the  
simplified notation $L_i(-)$ in this paper as the groups we consider have 
the required properties.

We conclude the introduction by mentioning that in fact we prove the 
Fibered Isomorphism conjecture in L$^{\l -\infty \r}$-theory for 
a more general class of groups, namely for any finite extension 
$\G$ of a {\it strongly poly-free group} (
see [\cite{AFR}, definition 1.1] or Definition 2.1 below) and deduce the 
isomorphism of the classical assembly map in $L$-theory for any  
torsion-free subgroup of $\G$. As a consequence we compute the surgery 
groups of the 
fundamental group of any fiber-type hyperplane arrangement   
complement in the complex $n$-space ${\Bbb C}^n$ (see 
Theorem \ref{mainthm2}). 
 
\section{Statements of the Main Theorem and its consequences}

Let us recall the definition of the strongly poly-free groups.

\begin{defn} (\cite{AFR}) {\rm A discrete group $\G$ is called 
{\it strongly poly-free} if there exists a finite filtration of $\G$ 
by subgroups: 
$1=\G_0\subset \G_1\subset \cdots \subset \G_n=\G$ 
such that the following conditions are satisfied:

1. $\G_i$ is normal in $\G$ for each $i$

2. $\G_{i+1}/\G_i$ is a finitely generated free group 

3. for each $\g\in \G$ and $i$ there is a compact 
surface $F$ and a diffeomorphism $f:F\to F$ such that the induced 
homomorphism $f_{\#}$ on $\pi_1(F)$ is equal to $c_\g$ in 
$Out(\pi_1(F))$, where $c_\g$ is the action of $\g$ on 
$\G_{i+1}/\G_i$ by conjugation and $\pi_1(F)$ is
identified with $\G_{i+1}/\G_i$ via a suitable isomorphism.

In such a situation we say that the group $\G$ has {\it rank} $\leq
n$.} \end{defn}

We now state our main theorem. 

\begin{thm}\label{mainthm} Let $\G$ be a finite extension of a strongly poly-free 
group (the finite group is the quotient group). 
Then the Fibered Isomorphism conjecture of Farrell and Jones 
in $L^{\l -\infty \r}$-theory is true for any subgroup of $\G$. In particular, 
it is true for any subgroup of $B_n$. 
\end{thm}

Although we prove Theorem \ref{mainthm} for the conjecture in $L^{\l -\infty \r}$-theory 
stated in [\cite{FJ}, \S1.7], the proof goes through, under certain 
conditions (see [\cite{R1}, $3(b)$ of theorem 2.2]), in a general setup 
of the conjecture in equivariant homology theory formulated in 
\cite{BL0} and for a more general class of groups. 

A corollary of the above theorem and [\cite{FR}, theorem 1.1] is the following.
This shows that for any torsion-free subgroup $G$ of $\G$ the surgery group  
$L_i(G)$ is isomorphic to the generalized homology group  
$H_i(BG , {\Bbb L}_0)$. Here ${\Bbb L}_0$ is a $1$-connective $\Omega$-spectrum 
with $0$th space homotopy equivalent to the classifying space $G/TOP$. 
This spectrum 
and the {\it assembly map} (or {\it universal homomorphism}) mentioned 
in the below statement were originally 
constructed by Quinn (\cite{Q1}, \cite{Q2}) using geometric methods.  For an 
algebraic treatment on this subject see \cite{Ra1} and \cite{Ra2}.

\begin{cor}\label{main} The classical assembly map in surgery 
theory is an isomorphism for any torsion-free subgroup of $\G$. That is, 
$H_i(BG , {\Bbb L}_0)\to L_i(G)$ is an isomorphism for all $i$ and 
for all torsion-free subgroups $G$ of $\G$. In particular, the 
assembly map is an isomorphism for any subgroup of $B_n$.
\end{cor}

\begin{proof} Let $H$ be a torsion-free group so that the following are 
satisfied.

$(1)$. $Wh(H)=K_{-i}({\Bbb Z}H)=\tilde K_0({\Bbb Z}H)=0$ for all $i\geq 1$.

$(2)$. The Isomorphism conjecture in $L^{\l -\infty \r}$-theory is true for $H$.

Then it is a known fact that for all $i$, 
$H_i(BH , {\Bbb L}_0)\to L_i(H)$ is an isomorphism. For a 
detailed proof see [\cite{LR}, theorem 1.28] or [\cite{FJ}, 
1.6.3].

Now the proof of the Corollary is immediate since $(1)$ is satisfied for 
$G$ by [\cite{FR}, theorem 1.1] and $(2)$ is satisfied by 
Theorem \ref{mainthm}. 

The particular case follows since 
$PB_n$ is strongly poly-free and $B_n$ is torsion-free (see 
the discussion after the following Remark).\end{proof}

\begin{rem}{\rm Here we recall that the isomorphism 
of the above assembly map is expected when the group 
is torsion-free. The integral Novikov conjecture in 
$L$-theory states that this assembly map should 
be split injective.}\end{rem}

Before we state our main computation of the surgery groups
we recall the definition of a fiber-type 
hyperplane arrangement from [\cite{OT}, p. 162]. Such an arrangement 
${\cal A}_n\subset {\Bbb C}^n$, that is the union of a finite number of 
affine hyperplanes in ${\Bbb C}^n$ is called {\it strictly linearly fibered} 
if after a suitable linear change of coordinates, the restriction of the 
projection of ${\Bbb C}^n-{\cal A}_n$ to the first $(n-1)$ coordinates 
is a fiber bundle projection whose base space is the complement of an 
arrangement ${\cal A}_{n-1}$ in ${\Bbb C}^{n-1}$ and whose fiber is the 
complex plane minus finitely many points. By definition the arrangement 
$0$ in ${\Bbb C}$ is fiber-type and ${\cal A}_n$ is defined to be 
{\it fiber-type} if ${\cal A}_n$ is strictly linearly fibered 
and ${\cal A}_{n-1}$ is of fiber type. It follows by repeated application 
of the homotopy exact sequence of a fibration that the complement 
${\Bbb C}^n-{\cal A}_n$ is aspherical. And hence $\pi_1({\Bbb C}^n-{\cal A}_n)$ 
is torsion-free.

The hyperplane arrangement ${\cal H}_n$ for $PB_n$ 
as described in the Introduction is an example of a fiber-type 
arrangement.

Now recall from [\cite{FR}, theorem 5.3] that if 
${\cal A}$ is a fiber-type 
hyperplane arrangement in ${\Bbb C}^n$, then the fundamental group  
$\pi_1({\Bbb C}^n-\cup {\cal A})$ is strongly poly-free. In 
particular $PB_n$ is also 
strongly poly-free. This was proved in [\cite{AFR}, theorem 2.1]. 

As a consequence of Theorem \ref{mainthm} we prove the following.

\begin{thm}\label{mainthm2}
Let ${\cal A}=\{A_1, A_2,\ldots , A_N\}$ be a fiber-type 
hyperplane arrangement in ${\Bbb C}^n$, then the surgery groups 
of $\G=\pi_1({\Bbb C}^n-\cup_{j=1}^{N}{A_j})$ are  
given by the following.

$$
L_i(\G) =\begin{cases} {\Bbb Z}&\text{if
$i\equiv$ 0 mod 4}\\
{\Bbb Z}^N &\text{if $i\equiv$ 1 mod 4}\\

{\Bbb Z}_2&\text{if $i\equiv$ 2 mod 4}\\

{\Bbb Z}_2^N &\text{if $i\equiv$ 3 mod 4}.
\end{cases}$$

\end{thm}

\section{The Isomorphism Conjecture and related results}
The Isomorphism conjecture of Farrell and Jones ([\cite{FJ}, \S1.6, \S1.7]) 
is a fundamental conjecture and implies 
many well-known conjectures in algebra and topology (see \cite{Lu} for a quick 
introduction to the conjecture and its consequences or see \cite{LR}). The statement 
of the conjecture  has been stated in a general setup of equivariant 
homology theory in \cite{BL0}. We recall the statement below.

Let ${\cal H}^?_*$ be an equivariant homology theory with values in 
$R$-modules for $R$ a commutative associative ring with unit. 

A {\it family} of subgroups of a group $G$ is defined as a set   
of subgroups of $G$ which is closed under taking subgroups and 
conjugations. If $\cal C$ is a class of 
groups which is closed under isomorphisms and taking subgroups then we 
denote by ${\cal C}(G)$ the set of all subgroups of $G$ which belong to 
$\cal C$. 
Then ${\cal C}(G)$ is a family of subgroups of $G$. For example $\cal 
{VC}$, the class of virtually cyclic groups, is closed under isomorphisms 
and taking subgroups. By definition a {\it virtually cyclic group} has a 
cyclic subgroup of finite index. 
   
Given a group homomorphism $\phi:G\to H$ and a family $\cal C$ of 
subgroups of $H$ define $\phi^*{\cal C}$ to be the family 
of subgroups $\{K<G\ |\ \phi (K)\in {\cal C}\}$ of $G$. Given a family 
$\cal C$ of subgroups of a group $G$ there is a $G$-CW complex $E_{\cal 
C}(G)$ which is unique up to $G$-equivalence satisfying the property that 
for $H\in {\cal C}$ the fixpoint set $E_{\cal C}(G)^H$ is 
contractible and $E_{\cal C}(G)^H=\emptyset$ for $H$ not in ${\cal C}$. 

Let $G$ be a group 
and $\cal C$ be a family of subgroups of $G$. Then the {\it 
Isomorphism conjecture} for the pair $(G, {\cal C})$ states that the 
projection 
$p:E_{\cal C}(G)\to pt$ to the point $pt$ induces an isomorphism 
$${\cal H}^G_n(p):{\cal H}^G_n(E_{\cal C}(G))\simeq {\cal H}^G_n(pt)$$ for $n\in 
{\Bbb Z}$. 

And the {\it Fibered Isomorphism conjecture} for the pair $(G, {\cal 
C})$ states that for any group homomorphism $\phi: K\to G$ the 
Isomorphism conjecture is true for the pair $(K, \phi^*{\cal C})$.

In this article we are concerned with the equivariant homology 
theory arising in $L^{\l -\infty \r}$-theory and when ${\cal C}={\cal {VC}}$ 
and $R={\Bbb Z}$. 
This (Fibered) Isomorphism conjecture is  equivalent to the Farrell-Jones conjectures 
stated in ([\cite{FJ}, \S 1.7]) [\cite{FJ}, \S 1.6]. For details see 
[\cite{BL0}, \S5 and \S6].

We say that the FICwF$^L$ is true for a group $G$ 
if the Fibered 
Isomorphism conjecture in $L^{\l -\infty \r}$-theory is true for $G\wr H$ for 
any finite group $H$. Here $G\wr H$ denotes the semidirect 
product $G^H\wr H$ with respect 
to the regular action of $H$ on $G^H=G\t G\t \cdots \t G$ ($|H|$ number 
of factors). Also we say that the FIC$^L$ is true for a group $G$
if the Fibered
Isomorphism conjecture in $L^{\l -\infty \r}$-theory is true for $G$.

Next, we recall some standard results and some recent development 
in this area which we need for the proof of Theorem \ref{mainthm}. Also 
we prove some basic results. Let us start by recalling that the 
Fibered Isomorphism conjecture has {\it the hereditary property}, that is 
if it is true for a group then it is true for any of its 
subgroups.

\begin{lemma} \label{group} Let $G$ be a group acting properly 
discontinuously and cocompactly by isometries on a metric space $X$. 
Then for any finite group 
$H$ the group $G\wr H$ acts properly discontinuously and 
cocompactly by isometries on the product metric space  
$X^H=X\t X\t \cdots \t X$ ($|H|$ 
number of factors).\end{lemma}

\begin{proof} This follows from the proof of Serre's theorem in 
[\cite{B}, theorem 3.1, p.190-191].\end{proof}

An immediate corollary to the above Lemma is the following. 
Recall that a  CAT$(0)$-space is a connected simply 
connected metric space which is nonpositively curved in the  
sense of distance comparison. For example the universal cover $\tilde M$ of a 
closed nonpositively curved Riemannian manifold $M$ with respect 
to the lifted metric is CAT$(0)$. For some more 
information on this subject see \cite{BH}.

\begin{cor} \label{coro} If $G$ acts properly discontinuously and 
cocompactly on a CAT$(0)$-space, 
then for any finite group 
$H$, $G\wr H$ also acts  properly discontinuously and 
cocompactly on a CAT$(0)$-space.\end{cor}

\begin{proof} The proof is immediate as the product of two 
CAT$(0)$-spaces is again CAT$(0)$.\end{proof}

A group $G$ is called CAT$(0)$ if it 
acts properly and cocompactly by isometries 
on a CAT$(0)$-space. Hence if $M$ is as above then $\pi_1(M)$ 
is a CAT$(0)$-group.
Therefore, by Corollary \ref{coro} for any finite group $H$, 
$\pi_1(M)\wr H$ is also a CAT$(0)$-group.
 
\begin{lemma}\label{product} The FIC$^L$ is true for $V_1\t V_2$ 
for any two virtually cyclic groups $V_1$ and $V_2$.\end{lemma}

\begin{proof} Since the FIC$^L$ is true for any virtually 
cyclic group we can assume that both $V_1$ and $V_2$ are 
infinite. Hence $V_1\t V_2$ contains ${\Bbb Z}\t {\Bbb Z}$ (=$A$, say) 
as a finite index normal subgroup. By the algebraic lemma in \cite{FR} 
$V_1\t V_2$ is a subgroup of $A\wr H$, where $H=(V_1\t V_2)/A$. Let 
$T$ be a flat $2$-dimensional torus. Then by Corollary \ref{coro} 
$A\wr H$ is a CAT$(0)$-group.
Therefore by [\cite{BL}, theorem B] the FIC$^L$ is true for $A\wr H$ 
since the CAT$(0)$-space $\tilde T^H$ is finite dimensional. Here $\tilde T$ 
denotes the universal cover of $T$ with the lifted metric. Hence FIC$^L$ is true  
for $V_1\t V_2$ by the hereditary property.\end{proof}

\begin{lemma}\label{inverse}
Let $p:G\to Q$ 
be a surjective group homomorphism and assume that the 
FICwF$^L$ is true for $Q$, for ker$(p)$ and for $p^{-1}(C)$ for 
any 
infinite cyclic subgroup $C$ of $Q$. 
Then $G$ 
satisfies the FICwF$^L$.\end{lemma}

\begin{proof} The proof is immediate using Lemma \ref{product} 
and [\cite{R1}, lemma 3.4].\end{proof}

\begin{lemma}\label{free} Let $G$ be isomorphic to 
one of the following groups.

\begin{itemize}

\item The fundamental group of a closed nonpositively 
curved Riemannian manifold.

\item The fundamental group of a compact $3$-manifold $M$ with 
nonempty boundary so that there is a fiber bundle projection 
$M\to {\Bbb S}^1$.

\item A finitely generated virtually free group.
 
\end{itemize}

Then the FICwF$^L$ is true for $G$.\end{lemma}

\begin{proof} Since FIC$^L$ is true 
for all finite groups we can assume that $G$ is infinite.

Let $M$ be a closed nonpositively curved Riemannian manifold 
so that $\pi_1(M)\simeq G$. Then by Corollary \ref{coro} for 
any finite group $H$, $G\wr H$ is a CAT$(0)$-group and 
hence the FIC$^L$ is true for $G\wr H$ by [\cite{BL}, theorem B] 
since the CAT$(0)$-space $\tilde M^H$ is finite dimensional. This 
completes the proof of the first item.

Now we give the proof for the second item, 
then the third one will follow using the hereditary property 
of the Fibered Isomorphism conjecture. 
 
Let $S$ be a fiber of the 
fiber bundle $M\to {\Bbb S}^1$ with monodromy diffeomorphism 
$f:S\to S$. Then $M$ is diffeomorphic to the mapping torus 
of $f$. Therefore $M$ is a compact Haken $3$-manifold 
(that is an irreducible $3$-manifold which has a 
$\pi_1$-injective embedded surface, see \cite{He}) with 
boundary components of zero Euler characteristic. 
We now apply [\cite{L}, theorem 3.2 and 3.3] to get a 
complete nonpositively curved Riemannian metric in the 
interior of $N$ so that near the boundary the metric 
is the product flat metric, that is each end is isometric 
to $X\t [0, \infty)$ for some flat $2$-manifold $X$. Therefore 
if we take the double $D$ of $M$, we get a closed nonpositively 
curved Riemannian manifold. Hence by the first item FICwF$^L$ is 
true for $\pi_1(D)$ and consequently for $\pi_1(M)$ also 
by the hereditary property. 

This completes the proof of the Lemma.

\end{proof}

\section{Isomorphism of the Assembly map and computation of the surgery 
groups}

In this section we give the proofs of Theorems \ref{mainthm} and 
\ref{mainthm2}. 

\begin{proof}[Proof of Theorem \ref{mainthm}] Let $G$ be a finite 
index strongly poly-free normal subgroup of $\G$. Then by the 
algebraic lemma in \cite{FR} $\G$ can be embedded as a subgroup in 
$G\wr (\G/G)$. Therefore using the hereditary property it is 
enough to prove the FICwF$^L$ for any strongly poly free group. 
Hence we can assume that $\G$ is strongly poly-free.

The proof is by induction on the rank of $\G$ and 
the framework of the proof 
is same as that of the proof of 
[\cite{R1}, $3(b)$ of theorem 2.2]. 

If the rank is $\leq 1$ then 
$\G$ is a finitely generated free group and hence 
the theorem follows from the third item in Lemma \ref{free}.
 
Therefore assume that the rank of $\G \leq k$ and that the 
FICwF$^L$ is true for all strongly poly-free groups $\G$
 of rank $\leq k-1$.

Let $1=\G_0< \G_1 < \cdots < \G_k=\G$ be a filtration of 
$\G$. 

Consider the following exact sequence. $$1\to \G_1\to \G\to 
\G/\G_1\to 1.$$ Let $q:\G\to \G/\G_1$ be the above projection. 
The following assertions are easy to verify.
 
\begin{itemize}
\item $\Gamma/\Gamma_1$ is strongly poly-free and has rank $\leq 
k-1$. 
\item $q^{-1}(C)$ is a finitely generated free group or isomorphic to the 
fundamental group of a compact Haken $3$-manifold $M$ with nonempty boundary 
so that there is a fiber bundle projection $M\to {\Bbb S}^1$, 
where $C$ is either the trivial 
group or an infinite cyclic subgroup of 
$\G/\G_1$ respectively.
\end{itemize}

Now we can apply the induction hypothesis, Lemma \ref{inverse}, 
and Lemma \ref{free} to complete the proof of the Theorem. 

The particular case for $B_n$ follows as $PB_n$ is strongly 
poly-free ([\cite{AFR}, theorem 2.1]) and is an index $(n+1)!$ 
normal subgroup of $B_n$.

\end{proof}

To prove Theorem \ref{mainthm2} 
we need the following lemma regarding the topology of 
an arbitrary hyperplane 
arrangement complement in the complex $n$-space ${\Bbb C}^n$.

\begin{lemma}\label{suspension} The first suspension $\Sigma (
{\Bbb C}^n-\cup_{j=1}^N A_j)$ of the complement of a 
hyperplane arrangement ${\cal A}=\{A_1,A_2,\ldots , A_N\}$ 
in ${\Bbb C}^n$ is 
homotopically equivalent to the wedge of spheres 
$\vee_{j=1}^N S_j$ where $S_j$ is 
homeomorphic to the $2$-sphere ${\Bbb S}^2$ for $j=
1,2, \ldots , N$.\end{lemma}

\begin{proof} Let ${\cal A}=\{A_1,A_2,\ldots , A_N\}$ be an  
arrangement by linear subspaces of ${\Bbb C}^n$. Then in 
[\cite{S}, $(2)$ of proposition 8] it is proved that $\Sigma (
{\Bbb C}^n-\cup_{j=1}^N A_j)$ is homotopically 
equivalent to the following space. 
$$\Sigma (\vee_{p\in P}({\Bbb S}^{2n-d(p)-1}-\Delta P_{<p})$$

We recall the notations in the above display from [\cite{S}, p. 464].
$P$ is in bijection with ${\cal A}$ under a map, say $f$. Also 
$P$ is partially ordered by the rule that $p<q$ for $p,q\in P$ 
if $f(q)$ is a subspace of $f(p)$. $\Delta P$ 
is a simplicial complex whose vertex set is $P$ and chains 
in $P$ define simplices. $\Delta P_{<p}$ is the 
subcomplex of $\Delta P$ consisting of all $q\in P$ so that $
q<p$. $d(p)$ denotes the 
dimension of $f(p)$. The construction of the embedding of 
$\Delta P_{<p}$ in ${\Bbb S}^{2n-d(p)-1}$ 
is also a part of [\cite{S}, proposition 8].  

Now in the situation of the Lemma obviously 
$d(p)=2n-2$ and $\Delta P_{<p}=\emptyset$ for all $p\in P$. 
Finally note that the suspension of a wedge of spaces is homotopically 
equivalent to the wedge of the suspensions of the spaces. This completes 
the proof of the Lemma.\end{proof}

\begin{proof}[Proof of Theorem \ref{mainthm2}]
Let $X={\Bbb C}^n-\cup_{j=1}^NA_j$. Then by Lemma 
\ref{suspension} $\Sigma X$ is homotopically 
equivalent to $\vee_{j=1}^N S_j$, where $S_j$ is 
homeomorphic to the $2$-sphere ${\Bbb S}^2$ for $j=
1,2, \ldots , N$. Let $h_i$ be a generalized homology 
theory. Then $h_i(X)$ is computed as follows.
$$h_i(X)=h_i(*)\oplus h_{i-2}(*)^N.$$ Where 
$h_{i-2}(*)^N$ denotes the direct sum of 
$N$ copies of $h_{i-2}(*)$ and $`*'$ is a single point space.

Now replacing $h_i(-)$ by $H_i(-, {\Bbb L}_0)$ and 
using the following known computation we complete the 
proof of the theorem.
 
$$
H_i (*, {\Bbb L}_0) = \begin{cases} {\Bbb Z} & \text{if $i\equiv$ 0 mod 4}\\
0 & \text{if $i\equiv$ 1 mod 4}\\
{\Bbb Z}_2 & \text{if $i\equiv$ 2 mod 4}\\
0 & \text{if $i\equiv$ 3 mod 4}.\end{cases}
$$
\end{proof}

\bibliographystyle{plain}

\begin{thebibliography}{60}

\bibitem{AFR}
C.S. Aravinda, F.T. Farrell and S.K. Roushon,
\newblock Algebraic $K$-theory of pure braid groups, 
\newblock {\em Asian J. Math.},  4 (2000) 337-344.

\bibitem{EA}
E. Artin,
\newblock Theory of braids,
\newblock {\em Ann. of Math.} (2) 48 (1947), 101-126.

\bibitem{BL0}
A. Bartels and W. L\"{u}ck,
\newblock Isomorphism Conjecture for homotopy $K$-theory and groups acting 
on trees,
\newblock {\em J. Pure Appl. Algebra}, 205 (2006), 660-696.

\bibitem{BL}
\bysame,
\newblock The Borel conjecture for hyperbolic and 
$CAT(0)$-groups, preprint, arXiv:0901.0442v1[math.GT].

\bibitem{Bir}
J. Birman,
\newblock Braids, links and mapping class groups, Ann. of Math. Studies, 
\newblock Princeton University Press, Princeton, New Jersey, 1974.

\bibitem{BH}
M.R. Bridson and A. Haefliger,
\newblock Metric spaces of non-positive curvature,
\newblock Springer Verlag, Berlin, 1999. Die Grundlehren der 
mathematischen Wissenschaften, Band 319.

\bibitem{B} 
K. S. Brown, 
\newblock Cohomology of groups,
\newblock Springer-Verlag, New York Heidelberg Berlin, 1982. Graduate 
Texts in Mathematics 87.

\bibitem{FJ}
F.T. Farrell and L.E. Jones,
\newblock Isomorphism conjectures in algebraic $K$-theory, 
\newblock {\em J. Amer. Math. Soc.}, 6 (1993), 249-297.

\bibitem{FR}
F.T. Farrell and Sayed K. Roushon,
\newblock The Whitehead groups of braid groups vanish, 
\newblock {\em Internat. Math. Res. Notices}, no. 10 (2000), 515-526.

\bibitem{FN}
R. Fox and L. Neuwirth, 
\newblock The braid groups, 
\newblock {\it Math. Scand.} 10 (1962), 119-126.

\bibitem{He}
J. Hempel,
\newblock $3$-manifolds, 
\newblock {\em Annals of Mathematics Studies},
\newblock Princeton University Press, 1976.

\bibitem{KT}
C. Kassel and V. Turaev,
\newblock Braid groups. With the graphical assistance of Olivier Dodane. 
\newblock Graduate Texts in Mathematics, 247. Springer, New York, 2008.

\bibitem{L}
B. Leeb, 
\newblock 3-manifolds with(out) metrics of nonpositive curvature,
\newblock {\em Invent. Math.} 122 (1995), 277-289.

\bibitem{Lu}
W. L\"{u}ck,
\newblock The Farrell-Jones conjecture in algebraic $K$- and $L$-theory,
\newblock {\em Enseign. Math.} (2) 54 (2008), 140-141. In `Guido's book of conjectures',
a gift to Guido Mislin on the occasion of his retirement from ETHZ, June 2006. Collected by Indira Chatterji.
{\em Enseign. Math.} (2) 54 (2008), no. 1-2, 3--189.

\bibitem{LR}
W. L\"{u}ck and H. Reich, 
\newblock The Baum-Connes and the Farrell-Jones conjectures in $K$- and
$L$-theory, 
\newblock  In {\em Handbook of K-theory} Volume 2, edited by E.M. Friedlander, D.R.
Grayson, 703-842, Springer, 2005.

\bibitem{OT}
P. Orlik and H. Terao,
\newblock Arrangements of Hyperplanes,
\newblock Springer-Verlag, Berlin, New York, 1992.

\bibitem{Q1}
F. Quinn,
\newblock A geometric formulation of surgery, 
\newblock 1970 {\em Topology of Manifolds} (Proc. Inst., Univ. of Georgia, 
Athens, Ga., 1969) pp. 500–511 Markham, Chicago, Ill.

\bibitem{Q2}
\bysame, 
\newblock $B_{(TOP_n)^{\sim}}$ and the surgery obstruction,
\newblock {\em Bull. Amer. Math. Soc.} 77 (1971), 596-600.

\bibitem{Ra1}
A. A. Ranicki,
\newblock The total surgery obstruction,
\newblock Lecture Notes in Math. 763, Springer-Verlag, New York 1979, 275-316.

\bibitem{Ra2}
\bysame,
\newblock Algebraic $L$-theory and topological manifolds, 
\newblock Cambridge Tracts in Mathematics, 102. Cambridge 
University Press, Cambridge, 1992.

\bibitem{R1}
S.K. Roushon,
\newblock Algebraic $K$-theory of groups wreath product 
with finite groups, 
\newblock {\em Topology Appl.} 154 (2007), 1921-1930.

\bibitem{S}
Ch. Schaper,
\newblock Suspensions of affine arrangements.
\newblock {\em Math. Ann.} 309 (1997), no. 3, 463--473.

\end{thebibliography}
\ifx\undefined\bysame
\newcommand{\bysame}{\leavevmode\hbox to3em{\hrulefill}\,}
\fi

\medskip

\end{document}